\documentclass[11pt,a4paper]{amsart}
\begin{document}
\newtheorem{thm}{Theorem}[section]

\newtheorem{lem}[thm]{Lemma}
\newtheorem{cor}[thm]{Corollary}
\newtheorem{prop}[thm]{Proposition}
\newtheorem{conj}[thm]{Conjecture}

\theoremstyle{definition}
\newtheorem{defn}[thm]{Definition}
\newtheorem{rem}[thm]{Remark}
\newtheorem{constr}[thm]{Construction}
\newtheorem{nota}[thm]{Notation}
\newtheorem{example}[thm]{Example}
\def\Q{\mathbb Q}
\def\P{\mathbb P}
\def\p{\mathfrak p}
\def\C{\mathbb C}
\def\OO{\mathcal O}
\def\End{\mathop{End}}
\def\codim{{\hbox{\rm codim}\,}}
\def\dim{{\hbox{\rm dim}\,}}
\def\<{\langle}
\def\>{\rangle}
\def\phi{\varphi}
\def\symn{{\mathfrak S}_n}
\def\symnm{{\mathfrak S}_{n/m}}

\title[Orbifold cohomology for global quotients]
{Orbifold cohomology for global quotients}
\author{Barbara Fantechi}
\address{
Dipartimento di Matematica e Informatica\\
Universit\`a degli Studi di Udine\\
Via delle Scienze, 206\\
33100 UDINE}
\email{fantechi@dimi.uniud.it}

\author{Lothar G\"ottsche}
\address{International Center for Theoretical Physics
\\Strada Costiera 11
\\34100 Trieste, Italy}
\email{gottsche@ictp.trieste.it}
\begin{abstract}
Let $X$ be an orbifold which is a global quotient of a 
manifold $Y$ by a finite group $G$. We construct a 
noncommutative ring $H^*(Y,G)$ with a $G$ action 
such that $H^*(Y,G)^G$ is the orbifold cohomology ring of $X$ defined by 
Chen and Ruan. When $Y=S^n$, with $S$ a  surface with trivial canonical class
and $G={\mathfrak 
S}_n$, we prove that (a small modification of) the orbifold cohomology 
of $X$ is naturally isomorphic to the cohomology ring of the Hilbert 
scheme $S^{[n]}$, computed by Lehn and Sorger.  
\end{abstract}
\maketitle
\section*{Introduction}

If a finite group $G$ acts on a complex manifold $Y$, the quotient
$Y/G$ has a natural structure of smooth orbifold, which we denote by 
$[Y/G]$. Originating in physics
(see \cite{DHVW1}, \cite{DHVW2}, \cite{Z})
cohomological invariants, in particular orbifold Euler number and 
orbifold Hodge numbers
 of the orbifold $[Y/G]$ have been defined and studied,
with the idea that they should coincide with the invariants of a nice
(crepant) resolution of singularities. In \cite{B-B} it was shown 
that this is indeed true.
Recently  the cohomology ring for the resolution of symplectic
singularities was determined in Thm.~1.8(ii) of \cite{EG}.

One specially interesting case which was often used to test the 
physicists' predictions before a general proof was available
\cite{G}, \cite{GS}, \cite{HH}, is 
that of the Hilbert scheme of points on a surface; namely, $S$ is a 
complex surface, $Y=S^n$ and $G={\mathfrak S}_n$, the symmetric group 
on $n$ letters acting in the obvious way. A crepant resolution of the 
symmetric product is then provided by the Hilbert scheme $S^{[n]}$ (or, 
in case $S$ is not algebraic, by the corresponding Douady space).

Chen and Ruan \cite{CR1}, \cite{CR2} have introduced an orbifold 
cohomology ring for any orbifold, whose Hodge numbers coincide with the 
orbifold Hodge numbers. 
In the special case of the Hilbert scheme, 
they conjectured that the orbifold cohomology ring should be isomorphic 
to the cohomology ring of the Hilbert scheme if the surface $S$ 
has trivial canonical class.

Since the first version of \cite{CR1} appeared, 
the cohomology ring of $S^{[n]}$ has been explicitly 
computed by Lehn and Sorger \cite{LS2} for 
a surface $S$ with trivial $c_1(T_S)$
(building on previous results of Nakajima \cite{N1}, \cite{N2},
Grojnowski \cite{Gr},
Lehn \cite{L} and Li, Qin and Wang \cite{LQW}).

This paper started as an attempt to verify Chen and Ruan's 
conjecture using the result of Lehn and Sorger.
In section 3 we show that the conjecture is 
essentially true, i.e. by introducing suitable signs into the ring
structure of the orbifold cohomology ring, one gets 
a canonical ring isomorphism
respecting also the duality pairing. 

In order to prove this, we first (in section 1) introduce, for 
a manifold $Y$ with 
the action of a finite group $G$, a cohomology ring $H^*(Y,G)$. The 
construction is inspired by Chen and Ruan's definition of the orbifold
 cohomology. The ring 
$H^*(Y,G)$ is not commutative, on the other hand it is often simpler 
then the orbifold cohomology ring. 

The group $G$ acts naturally on the ring $H^*(Y,G)$, and the 
$G$ invariant subring is naturally isomorphic to the orbifold cohomology 
of the quotient, as proven in section 2; there we also make some 
considerations on the more general version of Chen and Ruan's conjecture 
on the relationship between orbifold cohomology and ordinary cohomology 
of a crepant resolution in the Hyperk\"ahler case.

In the fourth section we compute the orbifold cohomology of
Beauville's generalized Kummer varieties; we expect that also in this case
the orbifold cohomology ring (again up to an explicit sign change) 
will be isomorphic to  the ordinary cohomology
of the crepant resolution.

We collect in the appendix a few
elementary results about certain Galois covers of smooth and nodal 
curves.

This paper started from conversations with K.S.~Narain, G.~Thompson and
M.S.~Narasimhan.
We would like to thank M.~Lehn and C.~Sorger for
sharing and discussing with us the results of their paper \cite{LS2}.
A particular thank goes to M.S.~Narasimhan for very useful discussions, 
in particular about parabolic bundles and the relation to the shift 
in cohomology degree.

{
\def\SS{{\mathfrak S}}
We were informed that a result close to Theorem \ref{hilring}, 
namely the existence
of a ring isomorphism between $H^*_o([S^n/\SS_n])[nd]$ and the $\SS_n$
invariant part of $H^*(S)[n]\{\SS_n\}$, has been obtained independently by
Uribe \cite{U}.
}

This paper is dedicated to Susanna \cite{PB}.

\section{Orbifold cohomology for a quotient orbifold}

\subsection*{Notation and conventions}
A symbol $:=$ means that the left hand side is defined by the right hand
side.\par
All manifolds will be complex (although one could consider almost
complex ones with no major changes). All vector bundles will be complex
holomorphic. By dimension (or codimension) we will always mean complex
dimension (or codimension).
\par
For elements $g_1,\ldots,g_n$ in a
group $G$ we write $\<g_1,\ldots,g_n\>$ for the subgroup they generate.
\par
All group actions will be left group actions.
For a manifold $T$ with the action of a finite group $H$,
denote by $T^H$ the $H$-invariant locus (which is always a
closed submanifold); for every $g\in H$ let $T^g$
be the $g$-invariant locus. Write $T^{g_1,\ldots,g_n}$ for
$T^{\<g_1,\ldots,g_n\>}$.
\par
For a topological space $T$ let $H^{i}(T)$ denote
$H^{i}(T,\Q)$.
A morphism $f:T\to S$ of manifolds
induces via pullback
a degree-preserving ring homomorphism $f^*:H^{*}(S)\to H^{*}(T)$.
If $T$ is a submanifold of $S$ and $\alpha\in H^*(S)$, we write
$\alpha|_T$ instead of $i^*\alpha$, where $i:T\to S$ is the inclusion.
\par
If $f$ is proper, then
the pushforward  $f_*:H^{*}(T)\to H^{*}(S)$ on cohomology is
given as follows: to every cohomology
class $\alpha$, associate its Poincar\'e dual homology class $\bar\alpha$; if
$T$ is not
compact, this is a class in the Borel-Moore homology group. Then, as
$f$ is proper, $f_{*}$ is defined on Borel-Moore homology, and
the Poincar\'e dual of $f_*(\bar\alpha)$
in cohomology is defined to be
$f_{*}(\alpha)$.\par
Note that if $f:T\to S$ is an isomorphism, then
$f_*\alpha=(f^{-1})^*\alpha$.
We will denote the cup product of cohomology classes $\alpha$, $\beta$ just
by $\alpha\cdot\beta$ or $\alpha\beta$.
\par

\subsection*{The ambient ring as vector space with $G$-action}

Let $Y$ be a complex manifold with the action of a finite group $G$.
\begin{defn}
The vector space $H^*(Y,G)$ is defined as follows:
$$H^*(Y,G):=\bigoplus\limits_{g\in G}
H^{*}(Y^g).$$
\end{defn}

For $g\in G$ and $\alpha$ in $H^*(Y^g)$,
denote by $\alpha_g$ the corresponding element in the $g$-th
direct summand of $H^*(Y,G)$.
\begin{defn}\label{Gact}
If $g,h$ are two elements of $G$, then $h(Y^g)=Y^{hgh^{-1}}$.
Hence $G$ acts on $H^*(Y,G)$ by
$$h(\alpha_g):=(h_*\alpha)_{hgh^{-1}}.$$
\end{defn}

\begin{rem} There is an alternative way to define $H^*(Y,G)$ together with
the $G$ action. Let $(p,a):G\times Y\to Y\times Y$ be the map defined by
projection and action, and let $\bar Y=(p,a)^{-1}(\Delta_Y)$.
Then $(p,a)$ is $G$-equivariant with respect to the action
$h(g,y)=(hgh^{-1},hy)$ on $G\times Y$; as $\Delta_Y$ is $G$-invariant, so
is $\bar Y$.
The vector space $H^*(Y,G)$ can be viewed as
as $H^*(\bar Y)$, with the induced $G$-action. 
\end{rem}

\begin{rem}
The invariant subspace $H^*(Y,G)^G$ under the action of $G$ is isomorphic
to $$
\bigoplus\limits_{g\in T} H^*(Y^g)/C(g)$$
where $T\subset G$ is a set of representatives of the conjugacy classes of
$G$ (in particular $T=G$ if and only if $G$ is commutative) and $C(g)$ is
the centralizer of $g$ in $G$. Hence as a vector space $H^*(Y,G)$
coincides with
the orbifold cohomology of the quotient orbifold of $Y$ by $G$ as defined
in \cite{CR1}.\end{rem}

\subsection*{The grading}

The following definition \cite{Reid} 
is now becoming a standard, sometimes under
different names such as fermionic shift number \cite{Z} or degree shifting number
\cite{CR1}.

\begin{defn} Let $Y$ be a manifold of dimension $D$ with the action of a
finite group $G$. For $g\in G$ and $y\in Y^g$, let
$\lambda_1,\ldots,\lambda_D$ be the eigenvalues of the action of $g$ on
$T_{Y,y}$; note that they are roots of unity. Write $\lambda_j=e^{2\pi i
r_j}$ where $r_j$ is a rational number in the interval $[0,1[$. The {\em
age} of $g$ in $y$ is the rational number $a(g,y):=\sum_{j=1}^Dr_j$.
\end{defn}

The age $a(g,y)$ is a nonnegative rational number, and it is zero if
and only if $g$ acts as the identity in a neighborhood of $y$;
it is an integer if the action of $g$ near $Y$ preserves the
canonical class (i.e., the induced automorphism of $T_{Y,y}$ has
determinant $1$).

\begin{rem}\label{ageinv} The age $a(g,y)$  only depends on the connected
component $Z$ of $Y^g$ in which
$y$ lies; we can therefore denote it by $a(g,Z)$.
It is easy to check that
$$a(g,Z)+a(g^{-1},Z)=\codim(Z\subset
Y).$$\end{rem}

\begin{defn} We define a (rational) grading on $H^*(Y,G)$ as follows.
Let
$g\in G$ and let $Z$ be a connected component of $Y^g$, and $j:Z\to Y^g$
the inclusion. Let $\alpha\in H^i(Z)$; we assign to $j_*\alpha_g$ the
degree $i+2a(g,Z)$.
\end{defn}
Note that $H^*(Y,G)$ is integrally graded if the
age of every element of
$G$ at every point in its fixed locus is an integer (or in fact a
half-integer, as will be the case in the section 3).

\begin{defn}\label{evenandodd}
For later use, we also define a splitting of $H^*(Y,G)$ into even and odd
part, as follows: $$H^{ev}(Y,G)=\bigoplus_{g\in G}H^{ev}(Y^g)$$
and analogously for $Y^{odd}$.
\end{defn}

Note that $H^{ev}(Y,G)$ coincides with the
even-graded part if and only if for every $g\in G$ and for every $y\in
Y^g$ the age of $g$ in $y$ is an integer.

Note also that the $G$ action on $H^*(Y,G)$ preserves both the splitting
into even and odd parts and the grading.

\subsection*{Definition of the classes $c(g,h)$}

\par
We construct cohomology classes $c(g,h)\in H^*(Y^{g,h})$ which will be
used in defining the multiplication in $H^*(Y,G)$.

We use the following convention: a vector bundle on a disjoint
union of
manifolds
is the datum of a
vector bundle on each connected component, possibly having different ranks
on different components; its top Chern class is the cohomology class
restricting to the top Chern class on each connected component.

We also use the following fact: 
let $H$ be a finite group. 
Assume that $E$ is a bundle with $H$ action on a manifold $M$
on which $H$ acts trivially 
(i.e., we are given a homomorphism $H\to \End(E)$); we say that $E$ is an $H$-bundle.
Then 
the representation of $H$ 
defined by a fiber of $E$ is locally constant on $M$. 
In particular the $H$ invariant part of $E$ is also a vector bundle.

\begin{constr}\rm Let $g$ and $h$ be elements of $G$, and let
$H$ be the subgroup
they
generate; $H$ contains $(gh)^{-1}$.
Let  $C=C(\P^1,g,h,(gh)^{-1},H)$ be the induced 
Galois cover of $\P^1$, with Galois group
$H$, branched over the three
points $0,1,\infty$ (see the definition in the appendix). \par
The bundle $T_Y$ is $G$-equivariant over $Y$, hence 
$T_Y|_{Y^H}$ is an $H$ bundle.

Let $E$ on $Y^H\times C$ be the $H$-equivariant pullback
$E=\pi^{*}T_{Y}|_{Y^H}$, where $\pi$ is the projection of $Y^H\times C$ to
$Y^H$ (which is $H$-equivariant): this means that an element $f\in H$
maps the triple $(y,c,v)\in U\times C\times
E_{(y,c)}$ to $(y,f(c),f(v))$.
\par
\label{excessbdl} Define
$F(Y,g,h)$ to be $R^1\pi_*^H(E)$, where $R^{1}\pi_*^H$ is the derived
functor of the invariant pushforward. 
We will write just $F(g,h)$ when $Y$
is clear from the context. 
We define $c(g,h)$ to be the top Chern class of $F(g,h)$. 
Note that $F(g,h)$ is the $H$ invariant part of the $H$-bundle 
$R^1\pi_*(E)=T_Y|_{Y^H}\otimes R^1\pi_*(\OO_{C\times M})$, hence it is a 
vector bundle on $Y^H$.
\end{constr}

\begin{lem}\label{trivlem} (1) Let $g_1,g_2\in G$, and let 
$H=\<g_1,g_2\>$ and $g_3=(g_1g_2)^{-1}$.
The bundle $F(g_1,g_2)$ is isomorphic to 
$F(g_2,g_1)$ and to
$F(g_2,g_3)$. \par (2) Let $v\in G$, then $v:Y^{g,h}\to Y^{g',h'}$ 
where a prime denotes conjugation by $v$ (e.g., $g'=vgv^{-1}$). Then 
$v^*F(g',h')=F(g,h)$.
\end{lem}
\begin{proof} (1) The curve $C(\P^1,g_1,g_2,g_3,H)$ is isomorphic, as a 
curve with $H$ action, to $C(\P^1,g_{\sigma(1)}, 
g_{\sigma(2)},g_{\sigma(3)},H)$ for any permutation $\sigma$ of the 
indices $1,2,3$.\par (2) Conjugation by $v$ defines an 
isomorphism $\phi$ between $H$ and $H'=\<g',h'\>$. There is clearly 
a natural $\phi$ equivariant isomorphism between the curves
$C(\P^1,g,h,(gh)^{-1},H)$ and 
$C(\P^1,g,h,(gh)^{-1},H')$, and from this the result follows. 
\end{proof}
\begin{lem}\label{normbdl}
The same $F(Y,g,h)$ can be obtained by replacing in its definition
$E$ with $\tilde E=\pi^*N_{Y^H/Y}$.
\end{lem}
\begin{proof}
The subbundle $T_{Y^H}$ is the $H$-invariant part of $T_Y|_{Y^H}$; since
the group $G$ is finite,
it is a direct summand and $T_Y|_{Y^H}=T_{Y^H}\oplus N_{Y^H/Y}$ as
$H$-bundle.\par
On the other hand, $R^1\pi_*^H(\pi^*T_{Y^H})$ is zero, since
$H^1(C,\OO_C)^H$ is equal to $H^1(\P^1,\OO_{\P^1})$ which is zero.
\end{proof}

One of the fundamental facts that allows to define 
the orbifold cohomology ring, is that the
age $a(g)$ can be interpreted in terms  of parabolic bundles.

Let $u\in Y^H$. Then $E_u:=E|_{\{u\}\times C}$ is an $H$-equivariant bundle on $C$,
 giving
rise to a parabolic bundle on the quotient $\P^1=C/H$.
The formula of  \cite{Gro} for the parabolic degree $c^1(E_u,H)$ implies
\begin{equation}
c^1(E_u,H)=deg(p_*^H(E_u))+a(g,U)+a(h,U)+a((gh)^{-1},U),\label{pardeg}
\end{equation} 
where $p_*^H(E_u)$ is the $H$-equivariant pushforward via the natural 
projection $p:C\to \P^1$.

The following is the analogue of Lemma 4.2.2 in \cite{CR1}.

\begin{lem}\label{therank}
Let $U$ be a connected component of $Y^H$. Then
the coherent
sheaf $F(Y,g,h)|_U$ defined in Construction \ref{excessbdl}
is a vector bundle of rank 
\begin{equation}
a(g,U)+a(h,U)-a(gh,U)-\codim(U\subset Y^{gh}).\label{aeq}
\end{equation}
\end{lem}
\begin{proof}
Let $u\in U$. It is enough to show that the fiber $F_u$ of $F(Y,g,h)$ at
$u$  has dimension given by (\ref{aeq}). 
The $H$-equivariant bundle
$E_u:=E|_{\{u\}\times C}$ on $C$ is defined by a representation of $H$.
By the  principal theorem of \cite{Gro} it follows that the 
parabolic degree $c^1(E_u,H)$ is $0$.

Applying the Riemann-Roch theorem to (\ref{pardeg}) gives 
\begin{align*}\dim(\pi_*^H(E_u))-\dim(F_u)&=\chi(p_*^H(E_u))\\
&=\dim(Y)-a(g,U)-a(h,U)-a((gh)^{-1},U).
\end{align*}

The result follows by observing that $\pi_*^H(E_u)=T_{Y^H,u}$, and
by definition 
$a((gh)^{-1},U)=\dim(Y)-\dim(T_{Y^{gh},u})-a(gh,U)$.
\end{proof}

The following Lemma will be
needed in Section 3.
\begin{lem}\label{splitting}
Assume that there is a
product decomposition
$Y=Y_1\times Y_2$ such that
$H$ acts on each factor separately. Let $Y^H=Y^H_1\times Y^H_2$ be the
induced product decomposition: then
$F(g,h)=p_1^{*}F(Y_1,g,h)\oplus p_{2}^{*}
F(Y_2,g,h)$.
\end{lem}
\begin{proof}
The bundle $T_{Y}|_{Y^H}$ splits, as a vector bundle with $H$-action,
as direct sum $p_{1}^*T_{Y_1}|_{Y_1^H}\oplus p_{2}^*T_{Y_2}|_{Y^2_H}$.
Pulling back to $Y^H\times C$ and applying $R^{1}\pi_{*}^H$ both respect
this
direct sum decomposition; this proves the result.
\end{proof}

\subsection*{The product}

\begin{defn} Define a bilinear map $$
\mu:H^*(Y,G)\times H^*(Y,G)\to H^*(Y,G)$$
by $$
\mu(\alpha_g,\beta_h):=\gamma_{gh}$$
where $$\gamma=i_{*}(\alpha|_{Y^{g,h}}\cdot \beta|_{Y^{g,h}}\cdot
c(g,h))$$
and
$i:Y^{g,h}\to Y^{gh}$ is the natural inclusion.
\end{defn}

\begin{lem}
The bilinear map $\mu$ sends $H^i(Y,G)\otimes H^j(Y,G)$ to
$H^{i+j}(Y,G)$. Hence, it defines a graded multiplication (in general
non-commutative: its associativity will be proven in Theorem \ref{assoc})
on $H^*(Y,G)$.\end{lem}
\begin{proof}
Basically this follows from Lemma \ref{therank}.
$H^i(Y,G)$ is the direct sum, over $g\in G$ and over connected components
$Z$ of $Y^g$, of $H^{i-2a(g,Z)}(Z)$. It is therefore enough to prove the
following:
given two elements $g_1,g_2$ in $G$, connected components $Z_i$
in $Y^{g_i}$ and elements $\alpha_i\in H^{d_i}(Z_i)$, for every connected
component $U$ of $Z_1\cap Z_2$ the element
$$\iota_*(\alpha_1|_U\cdot\alpha_2|_U\cdot c(g_1,g_2)|_U)$$ is in
$H^{i+j-2a(g_1g_2,U)}(Y^{g_1g_2})$. Here $\iota:U\to Y^{g_1g_2}$ is the
natural inclusion.
As restriction preserves the grading, and $\iota_*$ raises the degree by
$2\,\codim(U\subset Y^{g_1g_2})$, the only thing left to prove
is indeed
Lemma \ref{therank}.\end{proof}

\begin{rem}
The map $\mu$ is equivariant with respect to the action of $G$.
In fact, this follows directly from Lemma \ref{trivlem} (2).
The multiplication $\mu$ is commutative if the group $G$ is abelian, or 
more generally if every point of $Y$ has abelian stabilizer: this 
follows from Lemma \ref{trivlem} (1). 
\end{rem}

We recall here a particularly simple special case of the excess
intersection formula. This is certainly well known in many contexts: it 
is proven for instance in \cite{Q}, Proposition 3.3.

Let $S$ be a manifold, $S_1$ and $S_2$ closed submanifolds, and assume
that $U:=S_1\cap S_2$ is also a submanifold of $S$; let $j_i:S_i\to S$ and
$\iota_i:U\to S_i$ be the natural inclusions. The excess bundle
$E(S,S_1,S_2)$ of $U$ as
intersection of $S_1$ and $S_2$ in $S$ is a bundle on $U$ ``measuring'' how much
the intersection of $S_1$ and $S_2$ fails to be transverse
along $U$.
It is defined as the cokernel of the natural map
$N_{U/S_1}\to N_{S_2/S}|_U$ or, equivalently, of
$N_{U/S_2}\to N_{S_1/S}|_U$. In particular it is equivalent, in the 
Grothendieck
group of vector bundles on $U$, to $T_S|_U+T_U-T_{S_1}|_U-T_{S_2}|_U$.
We denote the top Chern class of $E(S,S_1,S_2)$ by $e(S,S_1,S_2)$.

For any cohomology class $\alpha\in H^*(S_1)$, the following excess
intersection formula holds in the cohomology ring of $S_2$:
$$j_2^*j_{1*}(\alpha)=\iota_{2*}(e(S,S_1,S_2)\cdot
\iota_1^*(\alpha)).$$

\begin{lem} \label{condass}
A sufficient condition for the map $\mu$ to define an
associative product
on $H^*(Y,G)$ is that, for every ordered triple of elements
$(g_1,g_2,g_3)\in G$,
the
following relation hold in the cohomology ring of $W=Y^{g_1}\cap
Y^{g_2}\cap
Y^{g_3}$:
\begin{equation}c(g_1,g_2)|_W \cdot c(g_1g_2,g_3)|_W \cdot e_{12}=
c(g_1,g_2g_3)|_W\cdot c(g_2,g_3)|_W \cdot e_{23},\label{eq:condass}
\end{equation}
where $e_{12}=e(Y^{g_1g_2},Y^{g_1,g_2},Y^{g_1g_2,g_3})$ and
$e_{23}=e(Y^{g_2g_3},Y^{g_1,g_2g_3},Y^{g_2,g_3})$.
\end{lem}

\begin{proof}
Write $i:W\to Y^{g_1g_2g_3}$ for the natural inclusion.
Using the excess intersection formula,  it is a straightforward
computation to
check that
the product $(\alpha_{g_1}\cdot \beta_{g_2})\cdot \gamma_{g_3}$ is equal
to $i_*\lambda_{g_1g_2g_3}$ where $\lambda\in H^*(W)$ is $$\lambda=
\alpha|_W\cdot \beta|_W \cdot \gamma|_W
\cdot c(g_1,g_2)|_W  \cdot
c(g_1g_2,g_3)|_W\cdot e_{12}.$$
The only thing one needs to remember is that the assumption that all the
$c(g,h)$ have even degree allows one to move them around in the
product freely.
Analogously, $\alpha_{g_1}\cdot (\beta_{g_2}\cdot \gamma_{g_3})$
is equal to $i_*\lambda'_{g_1g_2g_3}$ where $\lambda'\in H^*(W)$ is
$$\alpha|_W\cdot \beta|_W \cdot \gamma|_W
\cdot c(g_1,g_2g_3)|_W  \cdot
c(g_2,g_3)|_W\cdot e_{23}.$$
\end{proof}

\subsection*{Proof of associativity}

\begin{thm}\label{assoc}
The
bilinear map $\mu$ defines a graded, $G$-equivariant, associative
multiplication on $H^*(Y,G)$.
\end{thm}
\begin{proof}
The only thing left to prove is associativity: we will check that the
sufficient condition of Lemma \ref{condass} is
verified. 
We will construct on $W$ two vector 
bundles $F_L$ and $F_R$ such that the left hand side of 
(\ref{eq:condass}) is
$c_{top}(F_L)$ (Lemma \ref{ass1})
and the right hand side of (\ref{eq:condass}) is
$c_{top}(F_R)$ (Remark \ref{ass2}); then
(Proposition \ref{ass3}) we will prove that $F_L$ and $F_R$ are isomorphic.

\end{proof}

\begin{nota}
From now on we will just write $\alpha\cdot\beta$ or $\alpha\beta$ 
instead of $\mu(\alpha,\beta)$ for the product in $H^*(Y,G)$.
\end{nota}

Write $g_4$ for the unique element of the group $G$ such that 
$g_1g_2g_3g_4$ is the identity. Note that $Y^{g_1g_2}=Y^{g_3g_4}$
and that $F(g_1g_2,g_3)$ is isomorphic 
to $F(g_3,g_4)$ by Lemma \ref{trivlem} (1).
Note also that $\<g_1g_2,g_3\>=\<g_3,g_4\>$.

\begin{lem}\label{ass0} Let $F_L$ be a bundle on $W$. Then its top 
Chern class is equal to the left hand side of (\ref{eq:condass}) if 
$F_L$ is equivalent, in the Grothendieck group of vector bundles on $W$,
to 
$$F(g_1,g_2)+F(g_3,g_4)+T_{Y^{g_1g_2}}-T_{Y^{g_1,g_2}}-T_{Y^{g_3,g_4}}+T
_W$$
where we suppress the $|_W$ from the notation. 
Analogously, a vector bundle $F_R$ has as top Chern class the right hand side
of (\ref{eq:condass}) 
if it is equivalent to
$$F(g_2,g_3)+F(g_4,g_1)+T_{Y^{g_2g_3}}-T_{Y^{g_2,g_3}}-T_{Y^{g_4,g_1}}+T
_W$$
\end{lem}
\begin{proof} The left hand side is the top Chern class of the bundle 
$F(g_1,g_2)\oplus F(g_3,g_4)\oplus E(Y^{g_1g_2}, Y^{g_1,g_2}, 
Y^{g_1g_2,g_3})$. As equivalent vector bundles in the Gro\-then\-dieck group 
have the same Chern classes the result follows.
\end{proof}
Note that the second formula in the Lemma above can be obtained from 
the first by a cyclic permutation of the indices $(1,2,3,4)$ to 
$(2,3,4,1)$.

Let $H$ be the finite subgroup of $G$ generated by $g_1,g_2,g_3$. Note 
that $T_Y|_W$ is a $H$-bundle.
We now construct a nodal curve $C$ with an $H$ action such that the 
bundle
$F_L=(T_Y|_W\otimes H^1(C,\OO_C))^H$ has the properties claimed in
Lemma \ref{ass0}. A similar construction, replacing $(g_1,g_2,g_3)$ by 
$(g_2,g_3,g_4)$ will yield a vector bundle $F_R$ as in the Lemma.
\begin{constr}
Let $D$ be the union of two smooth rational curves $D'$ and $D''$ 
meeting transversally at a point $p$. Choose distinct marked points 
$p_1,p_2$ on $D'$ and $p_3,p_4$ on $D''$.
Let $C=C(D,g_i,H)$ be the associated
connected Galois $H$ cover of $D$ (see the appendix for a definition), 
branched over the $p_i$'s and the node $p$.
Let $C'$ be the inverse image of $D'$ and $C''$ 
the inverse image of $D''$.
Let $Z=C'\cap C''$. Write $C_{12}$ (resp.~$C_{34}$) for a connected 
component of $C'$ (resp.~$C''$) invariant under $\<g_1,g_2\>$ 
(resp.~$\<g_3,g_4\>$) and isomorphic 
to $C(D',g_1,g_2,(g_1g_2)^{-1},\<g_1,g_2\>)$ (resp.~$C(D'',
g_3,g_4,(g_3g_4)^{-1},\<g_3,g_4\>)$).
\end{constr}
\def\one{{\bf 1}}
\def\Ind#1#2{\mathop{Ind}{}^{#1}_{#2}}
We use the following notation. For a finite group $K$, we denote by 
$\one_K$ the trivial one dimensional representation of $K$; if $K'$ is 
a subgroup of $K$ and $U$ is a representation of $K'$, we denote by 
$\Ind{K}{K'}(U)$ the induced representation of $K$ (see e.g. 
\cite{FultonHarris} page 32).
\begin{lem} In the Grothendieck group of representations of $H$, the 
vector space $H^1(C,\OO_C)$ is equivalent to 
\begin{align*}
&\Ind{H}{\<g_1,g_2\>}H^1(C_{12},\OO_{C_{12}})+
\Ind{H}{\<g_3,g_4\>}H^1(C_{34},\OO_{C_{34}})+\one_H+\\
&\Ind{H}{\<g_1g_2\>}\one_{\<g_1g_2\>} -
\Ind{H}{\<g_1,g_2\>}\one_{\<g_1,g_2\>} 
-\Ind{H}{\<g_3,g_4\>}\one_{\<g_3,g_4\>}.
\end{align*}
\end{lem}
\begin{proof} Let $Z=C'\cap C''$. We have an exact sequence $$
0\to \OO_C\to \OO_{C'}\oplus \OO_{C''}\to \OO_Z\to 0;$$ 
as all sheaves 
involved are $H$-sheaves and the maps are $H$-equivariant, the long 
exact sequence of cohomology is an exact sequence of finite dimensional 
$H$ representations. As an $H$ representation, $H^0(C,\OO_C)$ is equal 
to $\one_H$, since $C$ is connected. As the connected component $C_{12}$ of $C'$ 
has stabilizer $\<g_1,g_2\>$, the space $H^0(C, \OO_{C'})$ is equal to 
$\Ind{H}{\<g_1,g_2\>}\one_{\<g_1,g_2\>}$ (see the appendix for details); 
analogously $H^0(C,\OO_{C''})$ 
is the representation 
$\Ind{H}{\<g_3,g_4\>}\one_{\<g_3,g_4\>}$ and $H^0(C,\OO_Z)$ is
$\Ind{H}{\<g_1g_2\>}\one_{\<g_1g_2\>}$. Moreover $H^1(C',\OO_{C'})$ is 
equal to $\Ind{H}{\<g_1,g_2\>}H^1(C_{12},\OO_{C_{12}})$
and analogously for $H^1(C'',\OO_{C''})$.
\end{proof}
\begin{cor}\label{ass1} The top Chern class of $F_L$ defined above 
is the left hand side of (\ref{eq:condass}) in Lemma \ref{condass}.
\end{cor}
\begin{proof} For any $H$ bundle $T$, any subgroup $K$ of $H$ and any 
representation $V$ of $K$, we have $(T\otimes \Ind{H}{K}V)^H=(T\otimes 
V)^K$. The result follows immediately by applying this remark to the 
induced representations in the previous Lemma, and taking $T=T_Y|_W$.
\end{proof}
\begin{rem}\label{ass2}
It is clear how to construct analogously  the bundle $F_R$: just 
permute cyclically everywhere in the definition of $F_L$  the indices 
$(1,2,3,4)$ to $(2,3,4,1)$, obtaining thus a curve $\bar C$ such that 
$F_R=(T_Y|_W\otimes H^1(\bar C,\OO_{\bar C}))^H$. It follows that the 
top Chern class of $F_R$ is the right hand side of (\ref{eq:condass}) in 
Lemma \ref{condass}.
\end{rem}
\begin{prop}\label{ass3} The bundles $F_L$ and $F_R$ are isomorphic. 
\end{prop}
\begin{proof} It is of course enough to prove that $H^1(C,\OO_C)$ is isomorphic as 
representation of $H$ to $H^1(\bar C,\OO_{\bar C})$. 
Consider the natural map $Ext^1(\Omega_C,\OO_C)^H\to {\mathcal 
E}xt^1(\Omega_C,\OO_C)$. It is easy to prove that it is surjective. As 
deformations of $C$ are unobstructed, this proves that there exists a 
smoothing of $C$ preserving the $H$ action; therefore there is a flat 
family of curves $f:\mathcal C\to B$ over a small disk such that $H$ acts 
fiberwise, the central fiber is $C$ and the other fibers are smooth. 
Hence the general fiber $\tilde C$ of $\mathcal C$ is isomorphic to 
the Galois $H$ cover $C(\P^1,g_1,g_2,g_3,g_4,H)$ of $\P^1$
branched over four points with stabilizers $g_1$, $g_2$, $g_3$ and 
$g_4$ (defined in the appendix); the $H$ representation 
$H^1(\tilde C, \OO_{\tilde C})$ does not 
depend on the particular branch points chosen. On the other hand, since 
$R^1f_*\OO_{\mathcal C}$ is an $H$ vector bundle with fiber the $H^1$ 
of the fiber, $H^1(\tilde C, \OO_{\tilde C})$ is isomorphic as a 
representation to $H^1(C,\OO_C)$. The same argument also applies to 
$\bar C$, completing the proof.

\end{proof}

\subsection*{The duality pairing}

\begin{defn}\label{duality}
Let $e\in G$ be the neutral element, assume that $Y$ is compact. Define
 a degree map $\int_{Y,G}:H^*(Y,G)\to \Q$
by $\int_{Y,G}\alpha_e=\int_Y \alpha$ for $\alpha\in H^*(Y)$ and
$\int_{Y,G}\beta_g=0$ for $\beta\in H^*(Y^g)$ and $g\ne 0$.
Define a duality pairing on  $H^*(X,G)$ by
$$\<\alpha_g,\beta_h\>=\int_{Y,G}\alpha_g\cdot\beta_h.$$
\end{defn}

\begin{lem}\label{pairing}
\begin{enumerate}
\item $\<\alpha_g,\beta_h\>=0$ if $h\ne g^{-1}$ and 
$\<\alpha_g,\beta_{g^{-1}}\>=\int_{Y^g} \alpha\beta$,
is just the Poincar\'e duality pairing on the compact manifold
$Y^g$.
\item The pairing $\<\, ,\, \>$ on $H^*(Y,G)$ is nondegenerate.
\end{enumerate}
\end{lem}
\begin{proof}
By definition $\alpha_g\cdot\beta_h=\gamma_{gh}$ for suitable $\gamma$, so 
$\<\alpha_g,\beta_h\>=0$ if $h\ne g^{-1}$.
By Lemma \ref{therank} combined with Remark \ref{ageinv}
we get that $F(Y,g,g^{-1})$ has rank $0$ and therefore $c(g,g^{-1})=1$.
Therefore $\alpha_g\cdot\beta_{g^{-1}}=i_*(\alpha\beta)$, where
$i:Y^g\to Y$ is the inclusion and $\alpha\beta$ is the cup product on $Y^g$.
Therefore $\<\alpha_g,\beta_{g^{-1}}\>=\int_{Y^g} \alpha\beta$.
This shows (1). (2) follows immediately from (1).
\end{proof}

\subsection*{Orbifold cohomology}

\begin{defn} The {\em orbifold cohomology} of the orbifold
$[Y/G]$ is the graded ring
$$H^*_o([Y/G]):=H^*(Y,G)^G.$$
It is a rationally graded associative ring.
\end{defn}

\begin{thm}\label{commutes} The orbifold cohomology is skew commutative
with respect
to the decomposition in even and odd part introduced in Definition
\ref{evenandodd}.
\end{thm}
\begin{proof}
Let $\tilde g,h\in G$; choose
$\tilde \alpha\in H^{n}(Y^{\tilde g})$ and $\beta\in H^{m}(Y^{h})$.
Define  $$\gamma:=
\sum_{f\in G}f^{-1}(\tilde \alpha_{\tilde g});$$ 
such classes $\gamma$ generate 
$H^{ev}(Y,G)^G$ for $n$ even and $H^{odd}(Y,G)^G$ for $n$ odd. 
We will prove that $\gamma\cdot \beta_h=(-1)^{mn} \beta_h\cdot\gamma$; this 
proves the stronger statement that every element of $H^*(Y,G)^G$ skew
commutes with every element
of $H^*(Y,G)$.

It is enough to check that, for any fixed $f\in G$,
$$f^{-1}(\tilde \alpha_{\tilde 
g})\cdot\beta_h=(-1)^{mn}\beta_h\cdot v^{-1}(\tilde\alpha_{\tilde g})$$ for $v=hf$. Let
$g=f^{-1}\tilde gf$,  $\alpha=f^*(\tilde\alpha)\in H^n(Y^g)$, and $H=\<g,h\>$.
Using this new notation, we have to check that
$$\alpha_g\cdot \beta_h=(-1)^{mn}\beta_h(h^*\alpha)_{h^{-1}gh}.$$

Let $H=\<g,h\>=\<h,h^{-1}gh\>$.
By definition $$\alpha_{g}\cdot \beta_{h}=i_{*}(\alpha|_{Y^{H}}\cdot
\beta|_{Y^{H}}\cdot c(g,h))_{gh},$$ where $i:Y^H\to Y^{gh}$ is the natural
inclusion; analogously
$$\beta_h\cdot (h^*\alpha)_{h^{-1}gh}=
i_*(\beta|_{Y^H}\cdot (h^*\alpha)|_{Y^H}\cdot c(h,g))_{gh}$$
Note that $c(g,h)=c(h,g)$
in view of Lemma \ref{trivlem} part (2); on the other hand, 
$(h^*\alpha)|_{Y^H}=\alpha|_{Y^H}$, and this completes the proof.
\end{proof}

\begin{defn} In case $Y$ is compact, define a degree map on 
$H^*_o([X/G])=H^*(X,G)^G$ by 
$\int_{[X/G]}:=\frac{1}{|G|}\int_{X,G}$.\par
In the same assumption, 
define a duality pairing on $H^*_o([X/G])$ by letting, for 
$\alpha,\beta\in H^*_o([X/G])$: 
$$\<\alpha,\beta\>_{[X/G]}=\int_{[X/G]}\alpha \cdot\beta
=\frac{1}{|G|}\<\alpha,\beta\>.$$ 
\end{defn}

\begin{rem}
Remark \ref{pairing} immediately implies that $\<\ ,\ \>_{X/G}$ is a nondegenerate
pairing. It is also easy to see that this pairing coincides with the pairing
defined in \cite{CR1}.
\end{rem}

We include here a few final remarks on the ring $H^*(Y,G)$.

\begin{rem} Let $Y$ be a complex manifold with an action by a finite group
$G$. Let $H$ be a subgroup of $G$ with the induced action. 
Then
$H^*(Y,H)=\sum_{h\in H} H^*(Y^h)$ is a subvector space of 
$H^*(Y,G)=\sum_{g\in G} H^*(Y^g)$, and from the definition it follows
immediately that $H^*(Y,H)$ is also a subring of $H^*(Y,G)$.
\end{rem}

In particular the results in section 3 give a very explicit description 
of $H^*(S^n,H)$ where $S$ is a smooth manifold and $H$ is any subgroup 
of ${\mathfrak S}_n$.

\begin{rem} Let $Y$ and $Z$ be two manifolds with the action of the 
same group $G$, $\phi:Y\to Z$ an \'etale $G$ equivariant map. Then $\phi$ 
induces a natural, degree preserving ring homomorphism $H^*(Z,G)\to 
H^*(Y,G)$ which is functorial. 
\end{rem}

In fact, the same functoriality property (pullback exists for \'etale 
maps) is true for the orbifold cohomology of arbitrary orbifolds, while 
pullback under general morphisms seems difficult to define. 
Thus orbifold cohomology has properties in between those of ordinary 
cohomology and those of quantum chomology (being, indeed, the degree 
zero quantum cohomology). 

As another instance of the closeness of orbifold cohomology to usual 
cohomology, the definition given in this paper 
can be modified to yield an orbifold Chow ring and a corresponding 
noncommutative ring $A^*(Y,G)$ (see \cite{F}).

\section{Orbifold cohomology and crepant resolutions}

\subsection*{Comparison with Chen and Ruan's definition}

\def\Hcr{H_{CR}}
In their paper \cite{CR1} Chen and Ruan define orbifold cohomology for 
an almost complex orbifold. In case the orbifold $X$ is a global quotient 
of a complex manifold $Y$ by a finite group $G$, their definition 
goes as follows. 

We start by remarking that, if in $G$ we have a relation $h=vgv^{-1}$, 
the element $v$ defines an isomorphism $v:Y^g\to Y^h$, hence a ring 
isomorphism
$v^*:H^*(Y^h)\to H^*(Y^g)$. If $h$ is also equal to $ugu^{-1}$, then 
$v^*$ differs from $u^*$ by $z^*$, where $z=uv^{-1}$ commutes with 
$h$. Let $C(g)$ be the centralizer of $g$; then the induced isomorphism 
$v^*:H^*(Y^h)^{C(h)}\to H^*(Y^g)^{C(g)}$ only depends on $g$ and $h$, 
and not on the choice of $v$. We denote it by $\iota_{h,g}$.

Therefore it makes sense to define $H^*(Y^g)^{C(g)}$ for any $g$ in a 
given conjugacy class $[g]$: a different choice of $g$ leads to a 
canonically isomorphic ring.

Chen and Ruan's orbifold cohomology is, as a vector space, 
$$\Hcr^*(X)=\bigoplus_{[g]\in T} H^*(Y^g)^{C(g)},$$
where $T$ is the set of all conjugacy classes of $G$.

Define a linear map $\psi: \Hcr^*(X)\to H^*(Y,G)^G$ by sending 
$\alpha\in H^*(Y^g)^{C(g)}$ to $\sum_{h\in [g]} (\iota_{g,h}\alpha)_h$. It 
is easy to check that $\psi(\alpha)$ is indeed $G$ invariant and $\psi$ 
is a linear isomorphism; it is also grade preserving, since our 
definition of grading is the same as in \cite{CR2}.

The product is defined via the pairing, as follows: for any three 
conjugacy classes $[g_1], [g_2], [g_3]$ and elements $\gamma_i\in 
H^*(Y^{g_i})^{C(g_i)}$ the triple pairing is $$
\<\gamma_1,\gamma_2,\gamma_3\>:=\sum_{[h_1,h_2]\in S} 
\frac{1}{|C(h_1,h_2)|}
\int_{X^{h_1,h_2}} \gamma_1'\cdot\gamma_2'\cdot \gamma_3'\, c(h_1,h_2)$$
where $S$ is the set of pairs $(h_1,h_2)\in [g_1]\times [g_2]$ such that 
$h_3=(h_1h_2)^{-1}\in [g_3]$ modulo 
simultaneous conjugation,  $C(h_1,h_2)=C(h_1)\cap C(h_2)$ and
$\gamma_i'=(\iota_{g_i,h_i}(\gamma_i))|_{Y^{h_1,h_2}}$.
Here $c(h_1,h_2)$ is the same cohomology class which was introduced in 
section 1.

The pairing $\<\gamma_1,\gamma_2\>$ is defined to be 
$\<\gamma_1,\gamma_2,\gamma_3\>$ where $g_3$ is be the identity of $G$ 
and $\gamma_3$ 
the identity in $H^*(Y)$ (i.e., the fundamental class of $Y$). The 
datum of the pairing and of the triple pairing defines a 
product by requiring that 
$\<\gamma_1\gamma_2,\gamma_3\>=\<\gamma_1,\gamma_2,\gamma_3\>$.

We want to prove that the two ring definitions coincide. To do this, 
let us compute $\int_{[X/G]}\psi(\gamma_1)\psi(\gamma_2)\psi(\gamma_3)$ 
and prove it coincides with  $\<\gamma_1,\gamma_2,\gamma_3\>$.

By definition, 
$$\int_{[X/G]}\psi(\gamma_1)\psi(\gamma_2)\psi(\gamma_3)=
\frac{1}{|G|}\sum_{(h_i)\in B} \int_{Y^{h_1,h_2}} \tilde 
\gamma_1\tilde \gamma_2\tilde \gamma_3 c(h_1,h_2)$$
where $\tilde\gamma_i=\gamma_i'|_{Y^{h_1,h_2}}$ and 
$B=\big\{(h_1,h_2,h_3)\bigm| h_i\in [g_i], \,h_3=(h_1h_2)^{-1}\}$.

The set $S$ above is equal to $B/G$, where $G$ acts by simultaneous conjugation.
For any $(h_1,h_2,h_3)\in B$, its stabilizer in $G$ is $C(h_1,h_2)$. 
Moreover, conjugation by elements of $G$ doesn't affect the product or 
the integral, so the sum over $(h_i)\in B$ is the same as the sum over 
$(h_i)\in S$ if we multiply the result by $|G|/|C(h_1,h_2)|$. Hence the 
two formulas give the same result. 

When we started working on this paper, the multiplication in
the orbifold cohomology of Chen and Ruan
differed from ours by a numerical factor. In the new version of their
paper the multiplication is changed
by changing the definition of the integral.

\subsection*{The conjecture}
Let $S$ be a complex surface with trivial canonical class (in other words, 
complex symplectic or Hyperk\"ahler). Let $X$ be the orbifold quotient 
of $S^n$ by the obvious action of the symmetric group ${\mathfrak 
S}_n$. 
The first part of Conjecture 6.3 in \cite{CR2} states that $\Hcr^*(X)$ 
coincides with the cohomology of the Hilbert scheme of $n$ points on 
$S$. Since the appearance of the first version of \cite{CR2}, the latter 
ring has been computed by Lehn and Sorger. We will prove in section 3 
that the conjecture holds modulo a certain sign change in the definition of 
the product. 

We want to make here some remarks on the second part of Conjecture 6.3 in  
\cite{CR2}, namely that if $X$ is an orbifold, $Z$ a crepant resolution of 
(the singular space associated to) $X$, such that 
both $X$ and $Z$ carry a Hyperk\"ahler structure, then 
the orbifold cohomology ring of $X$ coincides with the cohomology of $Z$.

\begin{nota} Until the end of this section $Y$ will be a complex manifold with 
the action of a finite group $G$, $X=[Y/G]$ the quotient orbifold, and 
$Z$ a crepant resolution of singularities of $Y/G$.
We will write $G'$ for $G\setminus \{e\}$, where $e$ is the identity of $G$.
\end{nota}
\begin{rem} Both $H^*(Z)$ and $H^*_o(X)$ have naturally a subring 
isomorphic to $H^*(Y/G)$; for $H^*(Z)$ it is the pullback, and for the 
orbifold cohomology it is $H^*(Y)^G_e$. If $Y$ is compact, the pairings 
also coincide. When we discuss the existence of an isomorphism between 
$H^*(Z)$ and $H^*_o(X)$, we will always require that it induces the 
identity on $H^*(Y/G)$.
\end{rem}

We compute the rings $H^*(Y,G)$ and $H^*_o([Y/G])$ 
in the following special case. In the rest of this section,
$Y$ will be a complex surface with the faithful Gorenstein action of a finite group 
$G$. Then $Y/G$ has 
rational double points as singularities, and its minimal resolution $Z$ is 
crepant. If $Y$ is a torus or a $K3$, then $Y$, $[Y/G]$ and $Z$ are 
Hyperk\"ahler.

Rational double points come in two series $A_n$ and $D_n$, plus three 
exceptional kinds $E_6$, $E_7$ and $E_8$. We describe in detail the $A_n$ 
case. Let $y\in Y$ be a point whose stabilizer $G_y$ is cyclic 
of order $n+1$ and let  $g_y$ be a generator of the stabilizer acting (in 
local coordinates $s,t$ on $Y$) by $g_y(s,t)=(\omega s,\omega^{-1}t)$, 
where $\omega=e^{2\pi i/(n+1)}$. Let $E(y)\in H^0(y)$ be the natural 
generator. Clearly given $g,h\in G_y$, in $H^*(Y,G)$ we have 
$E(y)_g\cdot E(y)_{h}=0$ unless 
$gh=e$, the identity of 
$G$; $E(y)_g\cdot E(y)_{g^{-1}}=p_e$, where $p\in H^4(Y)$ is the class 
of a point. In particular if $Y$ is compact 
$\<E(y)_g,E(y)_{g^{-1}}\>=1$ and $\<E(y)_g,E(y)_h\>=0$ if $h\ne -1$.
The quotient $Y/G$ has an $A_{n}$ singularity 
at the image of $y$. 
 
Assume that $G$ is cyclic of order $n+1$ and that $Y^g$ is the same for 
every $g\in G'$. 
Then the orbifold cohomology of $[Y/G]$ is isomorphic to 
$$
H^*(Y)^G\oplus\bigoplus_{y\in Y^G,g\in G'}\Q E(y)_g.$$

The cohomology of the resolution $Z$ is canonically isomorphic 
to $$H^*(Y)^G\oplus \bigoplus_{g\in G',y\in Y^g}\Q 
F(y)_g,$$ for suitable classes $F(y)_g\in H^2(Z)$, corresponding to the 
$(-2)$ curves.

At first sight the isomorphism between $H^*(Z)$ and $H^*_o([Y/G])$ is 
easy to construct.
However, in the cohomology of $Z$ the pairing
$\<F(y)_g, F(y)_{h}\>$ is equal to $-2$ if $g=h$, to $1$ if $h=gk$ or 
$g^{-1}k$ and is zero otherwise; in $H^*_o([Y/G])$
the pairing $\< E(y)_g, E(y)_{h}\>$ is equal to $1/(n+1)$ if $h=h^{-1}$ and 
$0$ otherwise.

If $n=1$ we can make $H^*(Z)$ isomorphic to a modified ring structure 
on
$H^*_o([Y/G],*)$;  it is enough to map
$F(y)_g$ to $2E(y)_g$ and change $c(g,g^{-1})$ from $1$ to $(-1)$.
A similar change of sign is required in the case of the Hilbert scheme 
of $n$ points on a surface, as we shall see in the next section.

However, no such generalization can be found if there is a point where 
$n\ge 2$; in fact, in $H^*(Z)$
the intersection pairing on the subspace generated by 
the $F(y)_g$ is negative definite, 
while each $E(y)_g$ is isotropic (i.e., $\<E(y)_g,E(y)_g\>=0$).
Hence even changing each class $c(g,h)$ by a sign or more generally by
a rational scalar will yield a non-isomorphic $\Q$-algebra. Note that 
the isomorphism must map the vector subspace generated by the $F(y)_g$'s 
to the vector subspace generated by the $E(y)_g$'s since in either case 
it's the orthogonal to the subspace $H^*(Y/G)$.

The same argument applies even if $X$ is not a global quotient, but 
just the smooth orbifold associated to a Gorenstein surface with 
rational double points.
Of course there are many $K3$ surfaces with an $A_n$ singularity 
with $n>1$, thus providing some kind of counterexample to the more 
general form of the conjecture of Chen and Ruan. 

We give an elementary example which is also a global quotient. Let 
$E$ be the elliptic curve which is a Galois triple cover of $\P^1$ 
branched over 
$3$ points, and let $f$ be a generator of the Galois group. Let 
$Y=E\times E$ and $G$ the automorphism group of $Y$ generated by 
$(f,f^{-1})$. Then $Y$ is Hyperk\"ahler, $G$ respects the Hyperk\"ahler 
structure, and the quotient $Y/G$ is a $K3$ surface with $9$ 
singularities of type $A_2$.

\section{The case of the symmetric product}

\def\enn{\{1,\ldots,n\}}
In this section, fix a smooth complex manifold $S$ of dimension $d$, and
a positive integer $n$.
Let $Y=S^n$, and $\symn$ the group of permutations of the set $\enn$; $\symn$
acts on $Y$ by $\sigma(s)_i=s_{\sigma(i)}$.
We prove that the orbifold cohomology $H^*_{o}([Y/\symn])$ is naturally
isomorphic up to a degree shifting to $H^*(S)^{[n]}$ as defined in
\cite{LS2}. There $H^*(S)^{[n]}$ is also shown, up to additional signs to be isomorphic
to $H^*(S^{[n]})$ in  case $S$ is a projective complex surface with 
numerically trivial canonical bundle.  

\begin{nota}\rm The notation introduced in Section 1 remains valid.\par
For a finite set $I$, denote by $S^I$
the manifold whose set of points is the set of maps from $I$ to $S$; it is
isomorphic to $S^r$, where $r=|I|$, the cardinality of $Y$.
In particular we identify $S^n$ with $S^{\enn}$.\par
A set map $\phi:I\to J$ induces a morphism $\tilde\phi:S^J\to S^I$
which is injective if $\phi$ is surjective and conversely: in the first
case it
is the inclusion of a multidiagonal, in the second it is a projection on
some of
the factors. Denote by $\tilde\phi_*$ and $\tilde\phi^*$ the induced
maps on cohomology.\par
For a subgroup $H$ of $\symn$, let $O(H)$ be the set of orbits
of $H$ in $\enn$; write $O(g)$ for $O(\<g\>)$ and 
$O(g,h)$ for $O(\<g,h\>)$. For $g\in \symn$, let $l(g)$ be
the minimal number of traspositions whose product is $g$.
Note that $n-l(g)$ is the cardinality of $O(g)$.
\par
For a graded vector space $V^*$, denote by $V^*[a]$ the graded
vector space defined by $V^i[a]=V^{i+a}$.
\end{nota}

\begin{rem} Let $H$ be a subgroup of $\symn$, and $\phi:\enn\to O(H)$ the
natural surjection. Then the image of $\tilde\phi:S^{O(H)}\to S^n$ is the
fixed locus of $H$.
\end{rem} 
\begin{lem} \label{agelem}
Let $g\in \symn$, $y\in Y^g$. The age of $g$ at $y$ is equal to $a(g)=d\cdot
l(g)/2$.
\end{lem}

In particular the age is always
a half-integer, and is an integer if $d=\dim S$ is even (the case which
will
interest us most being $\dim S=2$). The ring $H^*(Y,\symn)$ is therefore
integrally graded.

\begin{cor}\label{frank}
Let $g,h\in \symn$, then the rank of 
$F(g,h)$ is
$$r=\frac{d}{2}(n+2|O(g,h)|-|O(g)|-|O(h)|-|O(gh)|).$$
\end{cor}
\begin{proof}
This follows immediately from Lemma \ref{therank} and Lemma \ref{agelem}.
\end{proof}

For a graded Frobenius algebra $A$ (as defined in \cite{LS2}), 
let $A\{\symn\}$ be the Frobenius algebra defined in \cite{LS2}.
Its definition can be extended to the case where $A=H^*(S)[d]$ and $S$ 
is a noncompact manifold of dimension $d$  although in this case
no duality is defined on $A$;
it is enough to replace $e$ by $c_d(T_S)$ and to use the natural 
pushforward map on cohomology whenever needed.
\begin{prop}\label{isomo}
There is a canonical isomorphism of graded vector spaces
with $\symn$ action
$$\lambda:H^*(Y,\symn)[nd]\to H^{*}(S)[d]\{\symn\}.$$
If $S$ is compact, then the duality structures also agree.
\end{prop}
\begin{proof}
Both vector spaces are defined as direct sums over the elements of
$\symn$, so it is enough to define $\lambda$ componentwise.
As already remarked, the fixed locus of $g$ on $Y$ is naturally
isomorphic to the submanifold $S^{O(g)}$ of $S^{n}$ induced by the
natural surjection $\{1,\ldots,n\}\to O(g)$. This isomorphism
determines $\lambda$. The grading is preserved, since in both cases
it is chosen so as to have the graded pieces distributed symmetrical
around zero: in $H^{*}(Y,\symn)[nd]$ because $a(g)$ is equal to half
the real codimension of $Y^{g}$ in $Y$, and in $H^{*}(S)[d]\{\symn\}$
because in \cite{LS2} the grading of each summand of $A\{\symn\}$
is centered around zero.
The morphism $\lambda$ is $\symn$-equivariant by comparing Definition
\ref{Gact} with paragraph 2.8 in \cite{LS2}.

We see that $\lambda$ also preserves the duality, by comparing
Definition \ref{duality} with Proposition 2.16 in \cite{LS2}.
\end{proof}
In order to compare the product structure one has to compute the bundle
$F(g,h)$ defined in Construction \ref{excessbdl}. We begin by doing
so in a special case.

\begin{lem} Assume that $g,h$ are two elements in $\symn$ such that
$\<g,h\>$
acts transitively on
$\{1,\ldots,n\}$; in other words,
$Y^{g,h}$ is the small diagonal $\Delta$, canonically isomorphic to $S$.
Then the bundle $F(g,h)$ is isomorphic to a direct sum of copies of
$T_\Delta$.
\end{lem}\label{smalldiag}
\begin{proof}
We use the notation of Construction \ref{excessbdl}.
Let $V$ be the representation of $H$ on $\C^{n}$ induced
by the natural action of $\symn$.
As an $H$-equivariant vector bundle, $T_{Y}|_{\Delta}$ is isomorphic to
the tensor product of $T_{\Delta}$ (with the trivial $H$-action)
and of $V$.
This
implies that $F(g,h)=R^1\pi_*^H(E)$
is isomorphic to $T_\Delta\otimes W$, where $W$ is 
$H^{1}(C,\OO_C\otimes V)^{H}$.
\end{proof}
\begin{cor}\label{smdcor}
In the assumptions of the Lemma, $c(g,h)$ only depends
on the rank $r$ of $F(g,h)$; it 
has value $1$ if $r=0$, $c_{d}(T_\Delta)$ if $r=d$,
and is zero
otherwise. 
\end{cor}

\begin{thm} \label{hilring}
The linear isomorphism $\lambda:H^*(Y,\symn)[nd]\to
H^{*}(S)[d]\{\symn\}$
defined in Proposition \ref{isomo}
is a ring isomorphism.
\end{thm}
\begin{proof}
Let $h,l\in{\mathfrak S}_n$, and write
$H:=\<h,l\>$.  We have to prove that for every $\alpha\in H^*(Y^h)$ and 
$\beta\in H^*(Y^l)$ 
$$\lambda(\alpha_h\cdot\beta_l)=\lambda(\alpha_h)\cdot\lambda(\beta_h).$$
The fixed locus of $H$ is naturally isomorphic to $S^{O(H)}$: for every 
$o\in O(H)$, let $p_o:S^{O(H)}\to S$ be the natural projection.
Comparing with the definition of the product in Proposition 2.13 of 
\cite{LS2}, we see that we have to prove that 
\begin{equation}\label{keyeq}
c(h,l)=\prod_{o\in O(H)}p_o^*
\,\left(c_{d}(T_S)^{g(h,l)(o)}\right)\end{equation}
where $g(h,l)(o)$ is the graph defect defined in 2.6 in \cite{LS2}.
Note that $c_d(T_S)$ coincides with $e$ as defined in \cite{LS2}.
Because of the splitting Lemma \ref{splitting}, 
$$F(h,l)=\bigoplus_{o\in O(h)}p_o^*F_o(h,l)$$ where $F_o(h,l)=F(S^o,h,l)$
viewed as a bundle on $S$ (canonically isomorphic to $(S^o)^H$).

Hence is it is enough to prove (\ref{keyeq}) in
the case where $H$ acts transitively on $\{1,\ldots,n\}$; write 
$g(h,l)$ for $g(h,l)(\{1,\ldots,n\})$ (the unique orbit of $H$).
 
In this case the rank $r$ of $F(h,l)$ is equal, by comparing Corollary 
\ref{frank} with the definition of the graph defect, to $d\cdot g(h,l)$.
This completes the proof, by Corollary \ref{smdcor}, since $c_d(T_S)^i=0$
for every $i>1$.
\end{proof}

Now let $S$ be a projective surface over $\C$ with 
torsion canonical class (i.e., of Kodaira dimension zero). 
In \cite{LS2} Lehn and Sorger show that  after
introducing some additional signs
$H^*(S)[2]^{[n]}:=(H^*(S)[2]\{\symn\})^{\symn}$
is naturally isomorphic to $H^*(S^{[n]})[2n]$ as a ring. 

\begin{defn}\label{modring}
For $g,h\in \symn$ let $\epsilon(g,h):=(l(g)+l(h)-l(gh))/2$.
Note that this always is an integer.
We define a modified ring structure on $H^*(S^n,\symn)$ by.
$$\alpha_g * \beta_h=(-1)^{\epsilon(g,h)}\alpha_g\cdot\beta_h.$$
By the obvious identity $\epsilon(g,h)+\epsilon(gh,k)
=(l(g)+l(h)+l(k)-l(ghk))/2$
this defines an associative product.

This defines an induced ring structure on $H^*_o([S^n/\symn])$,
which we denote by $H^*_{so}([S^n/\symn])$, and we
define the pairing on $H^*_{so}([S^n/\symn])$ by 
$$\<\alpha,\beta\>:=\int_{[S^n/\symn]}\alpha*\beta.$$
\end{defn}

Let $g\in \symn$, and let $N=|O(g)|$. We identify as above $(S^n)^g$  with 
$S^{O(g)}$. Choosing a numbering $\phi:\{1,\ldots,N\} \simeq O(g)$ 
gives an isomorphism $\tilde\phi_*:H^*(S^N)\to H^*((S^n)^g)$.
So any class in $H^*((S^n)^g)$ can be written as
$\tilde\phi_*(\alpha_1\otimes \ldots\otimes \alpha_N)$.
We write 
$$\overline \phi(\alpha_1\otimes \ldots\otimes \alpha_N)
:=\sum_{h\in \symn}h(\tilde\phi_*(\alpha_1\otimes \ldots\otimes 
\alpha_N)_g)\in H^*_o([S^n/\symn]).$$

We denote by  $p_{k}:H^*(S^{[*]})\to H^*(S^{[*+k]})$ the 
generators of the Heisenberg algebra action on the cohomology of the 
Hilbert schemes (as defined e.g. in \cite{N1}).

Then, using the results of \cite{LS2}, Theorem \ref{hilring} can be
reformulated as follows.

\begin{thm}\label{hilrin} Let $S$ be a complex projective surface with $K_S=0$
and let $\one\in H^0(S^0)$ be the identity. 
There is a canonical ring isomorphism 
$\Psi:H^*_{so}([S^n,\symn])\to H^*(S^{[n]})$, given by
$$\overline \phi(\alpha_1\otimes \ldots\otimes \alpha_N)\mapsto
p_{|\phi(1)|}(\alpha_{1})\ldots p_{|\phi(N)|}(\alpha_{N}){\bf 1},$$ 
which is compatible with the duality pairing.
\end{thm}

\begin{proof}
Let $A$ be the Frobenius algebra $H^*(S)[2]$ with the degree map
given by $-\int_S$. Replace also the degree map on $S^{[n]}$  by 
$-\int_{S^{[n]}}$. Then in \cite{LS2} a ring isomorphism 
$A^{[n]}\to H^*(S^{[n]})[2n]$
compatible with the duality pairing is constructed.

In the definition of the product
$\lambda(\alpha_g)\cdot\lambda(\beta_h)$ on $A\{\symn\}$ the degree is used 
in two places. 
First, for the 
definition of the pushforward $i_*(\alpha|_{(S^n)^{g,h}}\cdot
\beta|_{(S^n)^{g,h}})$ where $i:(S^n)^{g,h}\to (S^n)^{gh}$ is the inclusion. 
With the new definition of the degree the pushforward changes by a sign 
$(-1)^{|O(g,h)|-|O(gh)|}$. 

Second, in the definition of the class $e$ in 
Section 2.2 of \cite{LS2}. In follows from the definition that 
$e$ is changed from $c_2(T_S)$ to $-c_2(T_S)$.
Let $h,l\in\symn$. Then, for the products to be the same,
$c_{top}(F(h,l))$ in (\ref{keyeq}) has to be  replaced by 
$(-1)^{rk(F(h,l))/2}c_{top}(F(h,l))$.

Putting this together, we have to change  $\alpha_g\cdot \beta_h$
by a factor of  $(-1)^{b(g,h)}$,
where by Lemma \ref{frank}
\begin{align*}
b(g,h)&=|O(gh)|-|O(g,h)|+\frac{1}{2}(n+2|O(g,h)|-|O(g)|-|O(h)|-|O(gh)|)\\
&=
\frac{1}{2}(n-|O(g)|-|O(h)|+|O(gh)|)=\frac{1}{2}(l(g)+l(h)-l(gh)),
\end{align*}
i.e. $\cdot$ is replaced by $*$.

The explicit formula in terms of the Heisenberg operators follows from the
definition of $\Phi$ in \cite{LS2} directly before Prop. 2.11.
\end{proof}

Note that, by deformation invariance of the cohomology ring, the 
assumption $S$ projective in Theorem \ref{hilrin} can be replaced 
by $S$ compact, since every compact complex surface with torsion 
canonical class can be deformed to a projective one. 
In fact, Theorem \ref{hilrin} 
is also true for $S={\mathbb A}^2$ by \cite{LS1}.

If  $S$ is a not necessarily compact 
complex surface which is also an abelian group
(e.g. an abelian variety), then the structure of 
the orbifold cohomology $H^*([S^n/\symn])$ is particularly simple.

\begin{cor}
Let  $S$ be a smooth complex surface which is also an abelian group.
Then the ring structure on $H^*(S^n,\symn)$ is given by
$\alpha_g\cdot\beta_h=\gamma_{gh}$, where
\begin{equation}\gamma=\begin{cases}
i_*(\alpha|_{(S^n)^{g,h}}\cdot \beta|_{(S^n)^{g,h}}) 
&\hbox{if }|O(g)|+|O(h)|+|O(gh)|=2|O(g,h)|+n,\\
0 & \hbox{ otherwise.
}\end{cases}\label{hilpro}
\end{equation}
 Here
$i$ is the embedding of $ (S^n)^{g,h}$ into $(S^n)^{gh}$.
\end{cor}
\begin{proof}
As all the Chern classes of $T_S$ vanish, we get $c(g,h)\ne 0$
if an only if the rank of $F(g,h)$ is zero, in which case $c(g,h)=1$.
By Lemma \ref{frank}  the rank of
$F(g,h)$ is $n+2|O(g,h)|- |O(g)|-|O(h)|-|O(gh)|$.
\end{proof}

We want to generalize the definition of $H^*_{so}([Y/G])$ 
from the case $Y=S^n$ for $S$ a surface with $K_S=0$ and
$G=\symn$ to 
arbitrary complex symplectic actions of a finite group $G$
on a complex symplectic manifold $Y$. This is based on the 
fact that in the above case $l(g)=a(g)$.

\begin{defn} 
Let $Y$ be a complex manifold with an action of a finite group $G$.
Assume that  for any pair 
of elements 
$g,h\in G$ with $Y^{g,h}$ nonempty, 
$\epsilon(g,h):=(a(g)+a(h)-a(gh))/2$ is an integer.

Then we can define a new associate ring structure on $H^*(Y,G)$ by
$\alpha_g *\beta_h =(-1)^{\epsilon(g,h)}\alpha_g\cdot \beta_h$.
This is associative because obviously 
$\epsilon(g,h)+\epsilon(gh,k)=(a(g)+a(h)+a(k)-a(ghk))/2$.
This defines a new ring structure on $H^*_o([Y/G])$, which we denote
by $H^*_{so}([Y/G])$.
\end{defn}





\section{Generalized Kummer varieties}

Now we want to compute the orbifold cohomology for the orbifold quotients
$[S_0^n/\symn]$
whose resolutions are the higher order Kummer varieties $K(S)_{n-1}$ of 
Beauville \cite{Bea}. 
In analogy with Theorem \ref{hilrin}, 
we expect that there is a canonical isomorphism
from  $H^*_{so}([S_0^n/\symn])$ to $H^*(K(S)_{n-1})$.

Assume that $S$ is a complex surface which is an abelian group, whose
identity element we denote by $0$.
Let $\sigma:=\sigma_n:S^{(n)}\to S$ be the morphism that associates
to a $0$-cycle its sum in $S$.
Let $\omega:S^{[n]}\to S^{(n)}$ be the Hilbert-Chow morphism. 
Then we define $K(S)_{n-1}:=\omega^{-1}\sigma^{-1}(0)$.
The most important case is if $S$ is compact, i.e. a $2$-dimensional 
torus, when one gets
the higher order Kummer varieties introduced and studied by Beauville 
in \cite{Bea}.
The varieties $K(S)_{n-1}$ are smooth and   complex symplectic: 
The proof  in \cite{Bea} works for any $S$, not necessarily compact
(of course in that case $K(S)_{n-1}$ will also be not necessarily compact).
 
Write $S^n_0:=\big\{(a_1,\ldots,a_n)\in S^n \bigm| \sum a_i=0\big\}$.
Then $S^n_0$ is isomorphic to $S^{n-1}$ and the symmetric group 
$\symn$ acts on $S^n_0$ by permuting the factors. $K(S)_{n-1}$ is a crepant
resolution of the quotient $S^n_0/\symn$. We now compute the orbifold cohomology
ring $H_o^*([S^n_0/\symn])$.
We describe $H^*(S)\times H^*(S_0^{n-1},\symn)$.
We denote by $S[k]$ the set of $k$-division points of $S$.
For any subgroup $H$ of $\symn$ let 
$m(H):=gcd\big\{|o|\bigm | o\in O(H)\big\}$ be the greatest common divisor
of the numbers of elements of the orbits of $H$ and 
$m(g_1,\ldots,g_r):=m(\<g_1,\ldots,g_r\>)$.

\begin{prop}\label{kumprop}
\begin{enumerate}
\item There is a canonical $\symn$-equivariant isomorphism
$$H^*(S)\times H^*(S_0^{n-1},\symn)\simeq \bigoplus_{g\in \symn}
\bigoplus_{y\in S[m(g)]} H^*((S^n)^g).$$
We denote by $\alpha_{g,y}$ a class $\alpha\in H^*((S^n)^g)$ in the summand
corresponding to $(g,y)$. Then the action of $\symn$ on the right hand side
is given by $h(\alpha_{g,y})=(h_*(\alpha))_{hgh^{-1},y}$.

\item 
The ring structure on $H^*(S)\times H^*(S_0^{n-1},\symn)$
induces via the above isomorphism the following ring structure
\begin{equation}\alpha_{g,x}\cdot \beta_{h,y}=\sum_{z\in S[m(g,h)]}
n_{g,h}(x,y,z) \gamma_{gh,z}.\label{kumpr}
\end{equation}
Here $\gamma\in H^*((S^n)^{gh})$ is given by (\ref{hilpro}) and 
$$\hbox{$n_{g,h}(x,y,z)=
\big| \big\{w\in S[m(g,h)]\bigm| \frac{m(g,h)}{m(g)}w=x,  
\frac{m(g,h)}{m(h)}w=y,\frac{m(g,h)}{m(gh)}w=z
\big\}\big|.$}$$
 
\end{enumerate}
\end{prop}
\begin{proof}
The proof will occupy the rest of this section.
We use some of the ideas of \cite{GS} p.~243. 
\begin{lem} Let $H$ be a subgroup of $\symn$ and assume $m(H)=1$. 
Then $(S^n_0)^H$ is connected and there is a canonical isomorphism  
$H^*(S\times (S^n_0)^H)\simeq 
H^*((S^n)^H)$.
\end{lem}
\begin{proof} 
There is an $\symn$-equivariant morphism $q:S\times S^n_0\to S^n$
 given on points
by $(a,(b_i)_i)\mapsto (a+b_i)_i$; for any subgroup $H$ of $ \symn$ 
its restriction is a morphism
$q:S\times (S^n_0)^H\to (S^n)^H$. The action of $S[n]$ on $S\times S^n_0$
by $c(a,(b_i)_i)=(a-c,(b_i+c)_i)$ commutes with the $\symn$-action
and the map $q$ is just the quotient map for this action.

Let $H\subset  \symn$ with $m(H)=1$. 
As in \cite{GS} p.~243 one shows that $S\times (S^n_0)^H$
 is isomorphic to $(S^n)^H$ (in particular $(S^n_0)^H$ is connected)
 and that the action of $S[n]$ on 
 $H^*(S\times (S^n_0)^H)$ is trivial. 
 Therefore $q^*$ is a  natural isomorphism $H^*((S^n)^H)\to 
 H^*(S\times (S^n_0)^H)$.
 \end{proof}

Let $H$ be a subgroup of $\symn$, and let $m=m(H)$.
 We identify $(S^n)^H$ with $S^{O(H)}$.
Note that the restriction of $\sigma$ to $(S^n)^H$ is given by
sending $(b_o)_{o\in O(H)}$ to $\sum_{o\in O(H)} |o|b_o$.
As all $|o|$ are divisible by $m$, we can define  
$$\sigma/m:(S^n)^H\to S,(b_o)_{o\in O(H)}\mapsto \sum_{o\in O(H)} 
\frac{|o|}{m}b_o,$$ and for $y\in S[m]$ we define
$(S^n)^H_y:=(\sigma/m)^{-1} y$.
By definition we get 
\begin{rem}
$(S^n_0)^H=\coprod_{y\in S[m]}(S^n)^H_y$.
\end{rem}
The proof of the next lemma shows  that this is the decomposition of $(S^n_0)^H$
into connected components.

\begin{lem}
There is a canonical isomorphism $H^*(S\times (S^n)^H_y)\simeq 
H^*((S^{n})^{H})$.
\end{lem}
\begin{proof}
 $(S^{n})^H_y$ is identified with 
 $$S^{O(H)}_y:=\big\{ (b_o)_{o\in O(H)}\in S^{O(H)}\bigm|
 \sum_{o\in O(H)} \frac{|o|}{m}=y\big\}.$$
 For any $z\in S$ with $\frac{n}{m}z=y$, we get an isomorphism 
 $S^{O(H)}_y\to S^{O(H)}_0,
(b_o)_o\mapsto (b_{o}-z)_o$. Grouping the elements of 
$\{1,\ldots,n\}$ in sets of $m$ elements, each of which contained
in an orbit of $H$, defines a surjection
$\{1,\ldots,n\}\to \{1,\ldots,n/m\}$, giving isomorphisms
$S^{O(H)}\simeq S^{O(\overline H)}$ and
 $S^{O(H)}_0\simeq S^{O(\overline H)}_0=(S^{n/m}_0)^{\overline H}$
for a subgroup $\overline H$ of $\symnm$ with $m(\overline H)=1$.

 Furthermore with the same proof as in 
 \cite{GS} p.~243 the induced isomorphism $H^*( (S^{n/m}_0)^{\overline H})
\to H^*((S^{n})^H_y)$
 is independent of the choice of $z$. As $m(\overline H)=1$,
  there is a canonical
 isomorphism $$H^*(S\times (S^{n/m}_0)^{\overline H})\simeq
 H^*((S^{n/m})^{\overline H})\simeq H^*((S^{n})^{H}).$$
 \end{proof}
 
 Putting everything together we get a canonical isomorphism
 $$H^*(S)\times H^*(S_0^{n},\symn)\simeq \bigoplus_{g\in \symn}
\bigoplus_{y\in S[m(g)]} H^*((S^n)^g),$$ which by definition commutes with
the $\symn$-action.
This shows part 1.

 We note that the action of
$\symn$ on $S^n_0$ is just the restriction of the action on $S^n$.
Furthermore for any subgroup $H\subset\symn$, the normal bundle of 
$(S^n_0)^H$ in $(S^n_0)$ is just the restriction of the 
 the normal bundle of 
$(S^n)^H$ in $(S^n)$.
In particular the age $a(g)$ of an element $g\in \symn$ is the same for 
both actions, and the bundle $F(g,h)$ on $(S^n_0)^{g,h}$ is the restriction of 
the corresponding bundle on $(S^n)^{g,h}$(in view of Lemma 
\ref{normbdl}).
Therefore the ring structure on $H^*(S_0^n,\symn)$ is determined in the same 
way as for $H^*(S^n,\symn)$:
$H^*(S_0^n,\symn)=\bigoplus_{g\in \symn} H^*((S_0^n)^g)$ and if we write $\alpha_g$
 for a class $\alpha\in  H^*((S_0^n)^g)$ in the summand corresponding to $g$,
 then 
 $\alpha_g\cdot\beta_h=\gamma_h$ with $\gamma$ given by the formula
(\ref{hilpro}) with the restriction to $(S^n)^{g,h}$ replaced by that over
$ (S^n_0)^{g,h}$.  Here
$i$ is the embedding of $ (S^n_0)^{g,h}$ into $(S^n_0)^{gh}$.

We need to determine how the product is distributed over the connected 
components $(S^n)^H_y$.

We write $\alpha_{(g,y)}$ for a class $\alpha\in H^*((S^n)^H_y)$
and $\beta_{g}$ for a class $\beta\in H^*((S^n_0)^H)$.
Let $g,h\in \symn$, $H:=\<g,h\>$, let $x\in S[m(g)]$, $y\in S[m(h)]$.
Then  
$\alpha_{(g,x)}\cdot \beta_{(h,y)}=0$ if 
$|O(g)|+|O(h)|+|O(gh)|\ne 2|O(g,h)|+n$.
Otherwise
\begin{align}
\alpha_{(g,x)}\cdot \beta_{(h,y)}
&=i_*(\alpha|_{(S^{n}_0)^H}\beta|_{(S^{n}_0)^H})_{gh}\nonumber\\
&=\sum_{w} (i_w)_*(\alpha|_{(S^{n})_w^H}\beta|_{(S^{n})_w^H})_{gh}\nonumber\\
&=\sum_{z\in S[m(gh)]}\sum_{w} 
(i_{w,z})_*(\alpha|_{(S^{n})_w^H}\beta|_{(S^{n})_w^H})_{(gh,z)}.\label{kumpro}
\end{align}
Here in the second line the inner sum is over all  $w\in S[m(H)]$
such that $(S^{n})_w^H\subset (S^{n})_x^{g}$ and 
$(S^{n})_w^H\subset (S^{n})_y^{h}$. In the last row we require in addition that
$(S^{n})_w^H\subset (S^{n})_z^{gh}$. Note that this is equivalent to 
$\frac{m(g,h)}{m(g)}w=x$,  
$\frac{m(g,h)}{m(h)}w=y$,
$\frac{m(g,h)}{m(gh)}w=z$.
Finally $i:(S^{n}_0)^H\to (S^{n}_0)^{gh}$, $i_w:(S^{n})_w^H\to (S^{n}_0)^{gh}$,
$i_{w,z}:(S^{n})_w^H\to (S^{n})_z^{gh}$ are the inclusions. 

From the definitions it is also obvious that, if $(S^{n})_y^H\subset 
(S^{n})_z^G$ and thus also $(S^n)^H\subset (S^n)^G$, then the isomorphisms
$H^*(S\times (S^{n})_y^H)\simeq H^*((S^n)^H)$ and 
$H^*(S\times (S^{n})_z^G)\simeq H^*((S^n)^G)$ commute with  the pullbacks induced by the 
inclusions.  Therefore (\ref{kumpr}) follows from (\ref{kumpro}),
thus completing the proof of Proposition \ref{kumprop}.
\end{proof}

\section*{Appendix}

We collect here for the reader's convenience a 
few facts about Galois coverings of rational smooth and nodal curves 
which are of an elementary nature and probably well known but
for which we couldn't find a 
suitable reference.

Let $D$ be a smooth complex curve and $U$ a contractible 
neighborhood of a point $p$. Then the fundamental group of $U\setminus 
p$ is infinite cyclic and has a canonical generator $\gamma_p$, defined 
by the condition that the integral over $\gamma_p$ of $dz/z$ be $2\pi i$, 
where $z$ is a local coordinate at $p$. Note that we don't need to 
specify a  basepoint since the fundamental group is abelian. 
Moreover, $\gamma_p$ doesn't depend on the neighborhood $U$ chosen, in 
an obvious sense.

If $D$ is complete and rational, given $n$ points $p_1,\ldots, p_n$ and a basepoint 
$p_0$ we can choose $\gamma_i\in 
F=\pi_1(D\setminus\{p_1,\ldots,p_n\},p_0)$ 
such that $\gamma_1\cdot\ldots\cdot\gamma_n$ is the identity and 
$\gamma_i$ is the pushforward of $\gamma_{p_i}$ from any contractible 
neighborhood of $p_i$ that doesn't contain any other $p_j$.
In fact $F$ is the free group generated by any $n-1$ of the $\gamma_i$.
So given any group $G$ and any set of $n$ elements $g_i$ such that 
$g_1\cdot\ldots\cdot g_n$ is the identity, there is an induced 
homomorphism $F\to G$ sending $\gamma_i$ to $g_i$; this defines a 
Galois covering $C^o(D,g_i,G)$ of $D\setminus\{p_1,\ldots,p_n\}$ 
with Galois group $G$.

If in the above construction $G$ is finite, then one can uniquely complete 
$C^o(D,g_i,G)$ to a proper smooth Galois $G$ cover $C=C(D,g_i,G)$ of $D$ branched over 
$p_1,\ldots,p_n$. Note that the representation of $G$ defined by 
$H^a(C,\OO_C)$ (where $a$ is $0$ or $1$) does not depend on the choice of $p_0,\ldots,p_n$ or the 
$\gamma_i$'s, because any two curves obtained by different choices can
be deformed to each other.

Note that $C$ is connected if 
the $g_i$'s generate $G$. Otherwise, if $H$ is the subgroup generated by 
the $g_i$, one can find a connected component of $C(D,g_i,G)$ which is 
$H$ invariant and isomorphic to $C(D,g_i,H)$ as a curve with $H$ action; 
the connected components of $C(D,g_i,G)$ are then in bijection with the 
cosets of $H$ in $G$.
This implies that the cohomology group $H^a(C(D,g_i,G),\OO)$
is isomorphic, as a representation of $G$, to 
$\Ind{G}{H}H^a(C(D,g_i,H),\OO)$ (the induced representation, see page 
32 of \cite{FultonHarris}).

Assume now that $D$ is the union of two smooth, proper rational curves 
$D'$ and $D''$ meeting transversely at one point $p$. Choose
$p_1,\ldots,p_k$ on $D'$ and $p_{k+1},\ldots,p_n$ on $D''$, and choose 
again $g_i\in G$ a finite group such that $g_1\cdot\ldots\cdot g_n$ is the 
identity.
Then the previous construction defines a Galois cover $C'=C(D',
g_1,\ldots,g_k,g,G)$
(resp.~$C''=C(D'',g^{-1},g_{k+1},\ldots,g_n,G)$) of $D'$ (resp.~$D''$) 
with Galois group $G$, branched over 
$p_1,\ldots,p_k,p$ (respectively $p,p_{k+1},\ldots,p_n$); 
here 
$g=g_{k+1}\cdot\ldots\cdot g_n=(g_1\cdot\ldots\cdot g_k)^{-1}$. 

Let $Z'$ (resp.~$Z''$) be the inverse image of $p$ 
in $C'$ (resp.~$C''$); then there are points $q'$ in $Z'$ and $q''$ in 
$Z''$ such that their stabilizer is generated by $g$ and such that $g$ 
acts on $T_{C',q'}$ and $T_{C'',q''}$ with two roots of unity with 
product one. One can therefore  naturally identify $Z'$ with $Z''$ by 
identifying $gq'$ with $gq''$ for every $g\in G$, thus obtaining a 
nodal curve $C=C(D,g_i,G)$ which is a Galois cover of $D$ with Galois
group $G$, branched over the marked points and the node. 
It is easy to see that $C$ is connected if and 
only if 
the $g_i$ generate $G$.

Again the definition of $C$ depends on a number of choices but the 
representation $H^a(C,\OO_C)$ of $G$ only depends on the elements 
$g_1,\ldots,g_n$.

The construction of the cover of the nodal curve
is closely related to the notion of admissible cover introduced in 
\cite{HarMum} and coincides with (a special case of)
that of twisted stable $n$-pointed map into $BG$ in \cite{AbrVis}.


\begin{thebibliography}{LQW}

\bibitem[AV]{AbrVis} D. Abramovich, A. Vistoli, {\em Compactifying 
the space of stable maps}, preprint math.AG/9908167 v2.

\bibitem[BB]{B-B}
V.~Batyrev, L.~Borisov, {\em Mirror duality and string-theoretic Hodge
numbers}, Invent. Math. 126
(1996), no. 1, 183--203.




\bibitem[Be]{Bea} A.~Beauville, {\em Vari\'et\'es K\"ahleriennes dont la
  premi\`ere classe de Chern est nulle}, J. Differential Geom. 18 (1983),
   no. 4, 755--782.

\bibitem[CR1]{CR1} W.~Chen, Y.~Ruan, {\em A new cohomology theory for
orbifold},
math.AG/0004129.
\bibitem[CR2]{CR2} W.~Chen, Y.~Ruan, {\em Stringy geometry and topology of
orbifolds}, math.AG/0011149.

\bibitem[D*1]{DHVW1} L.~Dixon, J.~A.~Harvey, C.~Vafa, E.~Witten,
{\em Strings on orbifolds},  Nuclear Phys. B 261 (1985), no. 4,
678--686.
\bibitem[D*2]{DHVW2}
L.~Dixon, J.~A.~Harvey, C.~Vafa, E.~Witten,
{\em Strings on orbifolds. II}, Nuclear Phys. B 274 (1986), no. 2,
285--314.

\bibitem[EG]{EG} P.~Etingof, V.~Ginzburg, {\em Symplectic reflection
algebras, Calogero-Moser space, and deformed Harish-Chandra homomorphism},
                 preprint math.AG/0011114.

\bibitem[F]{F} B.~Fantechi, {\em The orbifold Chow ring}, in 
preparation.

\bibitem[FH]{FultonHarris} W.~Fulton, J.~Harris, 
{\em Representation Theory: a first course}, GTM {\bf 129}, 
Springer-Verlag, New York 1991.

\bibitem[G]{G}
L.~G\"ottsche, {\em The Betti numbers of the Hilbert scheme of points
on a smooth projective surface}, Math. Ann. 286 (1990), no. 1-3, 193--207.



\bibitem[GS]{GS}
L.~G\"ottsche, W.~Soergel, {\em Perverse sheaves and the cohomology of
Hilbert schemes of smooth algebraic surfaces}. Math. Ann. 296 (1993), no. 2, 
235--245. 

\bibitem[Gr]{Gr} I.~Grojnowski, {\em Instantons and affine algebras.
  I. The Hilbert scheme and vertex operators}, Math. Res. Lett. 3
(1996), no. 2, 275--291.


\bibitem[Gro]{Gro} 
A.~Grothendieck, {\em Sur le m\'emoire de Weil: g\'en\'eralisation des 
fonctions ab\'eliennes}.
S\'eminaire Bourbaki, Vol. 4, Exp. No. 141, 57--71, 
Soc. Math. France, Paris, 1995. 

\bibitem[HM]{HarMum} J. Harris, D. Mumford,
{\em On the Kodaira dimension of the moduli space of curves}. 
With an appendix by William Fulton. 
Invent. Math. 67 (1982), no. 1, 23--88. 


\bibitem[HH]{HH} F.~Hirzebruch,  T.~H\"ofer,
{\em On the Euler number of an orbifold},  Math. Ann. 286 (1990), no. 1-3,
255--260.

\bibitem[L]{L} M.~Lehn, {\em Chern classes of tautological sheaves on
Hilbert schemes of points on surfaces}, Invent.
Math. 136 (1999), no. 1, 157--207.

\bibitem[LS1]{LS1} M.~Lehn, C.~Sorger,
{\em Symmetric groups and the cup product on the cohomology of 
Hilbert schemes},
       math.AG/0009131.



\bibitem[LS2]{LS2} M.~Lehn, C.~Sorger, {\em The cup product of the Hilbert
scheme for $K3$ surfaces},
math.AG/0012166.

\bibitem[N1]{N1} H.~Nakajima,  {\em Heisenberg algebra and Hilbert
schemes of points on projective surfaces},  Ann. of Math. (2)
145 (1997), no. 2, 379--388.

\bibitem[N2]{N2} H.~Nakajima, {\em Lectures on Hilbert schemes of points
on surfaces},  University Lecture Series, 18.
American Mathematical Society, Providence, RI, 1999. xii+132 pp.

\bibitem[LQW]{LQW}  W.~Li, Z.~Qin, W.~Wang,
{\em Vertex algebras and the cohomology ring structure of Hilbert schemes
of points on surfaces}, math.AG/0009132.



\bibitem[PB]{PB} {\em Susanna, prima bisiaca del millennio}. Il Piccolo,
edizione di Monfalcone, 2 gennaio 2000.

\bibitem[Q]{Q} D.~Quillen, {\em Elementary proofs of some results of 
cobordism using Steenrod Operations}, 
Adv. Math. {\bf 7} (1971), 29--56.

\bibitem[R]{Reid} M. Reid, {\em La correspondance de McKay}, S\'eminaire 
Bourbaki, exp.~867, vol.~1999/2000.

\bibitem[U]{U} B. Uribe, Ph.D.\ Thesis, University of Wisconsin, in 
preparation.

\bibitem[V]{V} E.~Vasserot,
{\em Sur l'anneau de cohomologie du sch\'ema de Hilbert de ${\mathbb
C}^2$,} math.AG/0009127.


\bibitem[Z]{Z} E.~Zaslow
{\em Topological orbifold models and quantum cohomology rings}.
Comm. Math. Phys. 156 (1993), no. 2, 301--331. 
\end{thebibliography}
\end{document}